\documentclass[12pt,reqno]{amsart}

\usepackage{amssymb}

\usepackage{mathrsfs}

\DeclareMathAlphabet{\mathpzc}{OT1}{pzc}{m}{it}

\usepackage[all]{xy}

\entrymodifiers={+!!<0pt,\fontdimen22\textfont2>}

\def\BQ{\mathbb{Q}}

\def\sD{\mathsf{D}}

\def\adots{\mathinner{\mkern1mu\raise1.0pt\vbox{\kern7.0pt\hbox{.}}\mkern2mu\raise4.0pt\hbox{.}\mkern2mu\raise7.0pt\hbox{.}\mkern1mu}}
\def\amp{\operatorname{amp}}

\def\C{\operatorname{C}}

\def\D{\sD}
\def\Dc{\D^{\operatorname{c}}}
\def\deg{\operatorname{deg}}

\def\dim{\operatorname{dim}}

\def\Dm{\D^-}
\def\Dp{\D^+}

\def\H{\operatorname{H}}
\def\hd{\operatorname{hd}}

\def\Hom{\operatorname{Hom}}

\def\inf{\operatorname{inf}}

\def\LTensor{\stackrel{\operatorname{L}}{\otimes}}

\def\Mod{\mathsf{Mod}}

\def\opp{\operatorname{o}}

\def\pcd{\operatorname{pcd}}
\def\pd{\operatorname{pd}}
\def\prod{\operatorname{prod}}

\def\RHom{\operatorname{RHom}}

\def\sup{\operatorname{sup}}

\numberwithin{equation}{section}

\newtheorem{Lemma}{Lemma}[section]
\newtheorem{Theorem}[Lemma]{Theorem}
\newtheorem{Proposition}[Lemma]{Proposition}
\newtheorem{Corollary}[Lemma]{Corollary}

\theoremstyle{definition}

\newtheorem{Setup}[Lemma]{Setup}

\newtheorem{Construction}[Lemma]{Construction}
\newtheorem{Remark}[Lemma]{Remark}

\newtheorem{Example}[Lemma]{Example}

\newtheorem{Notation}[Lemma]{Notation}
\newtheorem{Question}[Lemma]{Question}

\begin{document}

\setlength{\parindent}{0pt}
\setlength{\parskip}{7pt}
%The default \baselineskip is close to 4.8mm
%\setlength{\baselineskip}{5.3mm}

\title[Cochain DG algebras]
{Homological properties of cochain Differential Graded algebras}

\author{Anders J.\ Frankild}

\author{Peter J\o rgensen}
\address{School of Mathematics and Statistics,
Newcastle University, Newcastle upon Tyne NE1 7RU,
United Kingdom}
\email{peter.jorgensen@ncl.ac.uk}
\urladdr{http://www.staff.ncl.ac.uk/peter.jorgensen}

%\author{Next author goes here}
%\address{Next author's postal address goes here}
%\email{Next author's mail address goes here}

%\thanks{Date: \today. A thank you would go here}


%16E10: Homological dimension
%16E45: Differential graded algebras and applications
%16G10: Representations of Artinian rings 
%16G60: Representation type (finite, tame, wild, etc.) 
%16G70: Auslander-Reiten sequences (almost split sequences) and
%       Auslander-Reiten quivers
%18E30: Derived categories, triangulated categories
%18E35: Localization of categories 
%18G99: Homological algebra: None of the above, but in this section 
%55P62: Rational homotopy theory

\begin{abstract} 
  
  Consider a local chain Differential Graded algebra, such as
  the singular chain complex of a pathwise connected
  topological group.

\medskip
\noindent
  In two previous papers, a number of homological results were proved
  for such an algebra: An Amplitude Inequality, an
  Aus\-lan\-der-Buchs\-baum Equality, and a Gap Theorem.  These were
  inspired by homological ring theory.

\medskip
\noindent
  By the so-called looking glass principle, one would expect that
  analogous results exist for simply connected cochain Differential
  Graded algebras, such as the singular cochain complex
  of a simply connected topological space.

\medskip
\noindent
  Indeed, this paper establishes such analogous results.

\end{abstract}

\maketitle

\setcounter{section}{-1}
\section{Introduction}
\label{sec:introduction}

This paper is a sequel of \cite{FJiic} and \cite{PJampdgm}, or, more
accurately, their mirror image.  The papers \cite{FJiic} and
\cite{PJampdgm} investigated the homological properties of local chain
Differential Graded algebras, such as the singular chain complex of a
pathwise connected topological group.  Several results modelled on
ring theory were proved: An Amplitude Inequality, an
Aus\-lan\-der-Buchs\-baum Equality, and a Gap Theorem for Bass
numbers.

In this paper, we shall do the same thing for simply connected cochain
Differential Graded algebras, such as the singular cochain complex of
a simply connected topological space.  The resulting cochain
Aus\-lan\-der-Buchs\-baum E\-qua\-li\-ty and Gap Theorem are new,
while a cochain Amplitude Inequality was stated already in
\cite[prop.\ 3.11]{Schmidt}; our proof works by different methods.
For introductions to the theory of Differential Graded (DG) algebras,
we refer the reader to \cite{AFH}, \cite{FHTbook}, or \cite{KellerDG}.

One of the motivations for \cite{FJiic} was that the Gap Theorem
answered affirmatively \cite[Question 3.10]{AF} by Avramov and Foxby
on the so-called Bass numbers of local chain DG algebras.  The present
Gap Theorem implies that the answer is also affirmative for simply
connected cochain DG algebras.  In fact, it shows that for these
algebras, Avramov and Foxby's conjectural bound on the gap length of
the Bass numbers can be sharpened by an amount of one, see Corollary
\ref{cor:Bass_gap} and Remark \ref{rmk:Bass_gap}.

As indicated, the move from local chain to simply connected cochain DG
algebras is a general phenomenon: The so-called ``looking glass
principle'' of \cite{AH} states that each result on local chain DG
algebras should have a ``mirror i\-ma\-ge'' for simply connected
cochain DG algebras.  However, proofs cannot be translated in a
mechanical way.  A local chain DG algebra sits in non-positive
cohomological degrees and a simply connected cochain DG algebra sits
in non-negative cohomological degrees.  Accordingly, the simplest
statement of the looking glass principle is that it interchanges
positive and negative degrees.

For instance, over local chain DG algebras, in \cite{FJiic} we used
the notion of $k$-projective dimension of a DG module given by
\[
  \pd M = \sup \{\, j \,|\, \H_j(k \LTensor_R M) \not= 0 \,\}.
\]
The looking glass principle tells us that over simply connected
cochain DG algebras, we must replace this by
\[
  \pcd M = \sup \{\, j \,|\, \H_{-j}(k \LTensor_R M) \not= 0 \,\}
         = \sup \{\, j \,|\, \H^j(k \LTensor_R M) \not= 0 \,\}
\]
which will be called the {\em projective codimension} of $M$.  Indeed,
while the Auslander-Buchsbaum Equality for local chain DG algebras is
a statement relating $k$-projective dimension to other homological
invariants, see \cite[thm.\ 2.3]{FJiic}, for simply connected cochain
DG algebras it will be a statement relating the projective {\em
  co}dimension to other homological invariants, see Theorem
\ref{thm:AB} and Corollary \ref{cor:AB}.

The methods of this paper are different from the ones of \cite{FJiic}
and \cite{PJampdgm}.  They build on DG module adaptations of the ideas
used by Serre to prove \cite[thm.\ 10, p.\ 217]{Serre} on the
connection between homology and homotopy groups of topological spaces.
The main point is a list of termwise inequalities of power series
given in Propositions \ref{pro:inequalities} and
\ref{pro:compact_inequalities}.  An advantage of the present approach
is that, whereas the methods of \cite{FJiic} and \cite{PJampdgm} fail
for unbounded DG algebras, as demonstrated for example in \cite[sec.\ 
4]{FJiic}, this paper is able to treat bounded and unbounded DG
algebras on the same footing.

The paper is organized as follows: Section \ref{sec:background} gives
notation and e\-le\-men\-ta\-ry properties for DG modules over simply
connected cochain DG algebras.  Sections \ref{sec:construction} and
\ref{sec:inequalities} set up the DG module adaptation of Serre's
ideas from the proof of \cite[thm.\ 10, p.\ 217]{Serre}.  Section
\ref{sec:identities} proves the cochain Amplitude Inequality, the
Auslander-Buchsbaum Equality, and the Gap Theorem for Bass numbers in
Corollaries \ref{cor:amp}, \ref{cor:AB}, and \ref{cor:Bass_gap}.
These results arise as special cases of the stronger statements
Theorems \ref{thm:amp}, \ref{thm:AB}, and \ref{thm:gap}.

Finally, Section \ref{sec:topology} applies the Auslander-Buchsbaum
Equality and the Gap Theorem to the singular cochain DG algebra of a
topological space.  The context will be a fibration of topological
spaces, and we recover in Theorem \ref{thm:additivity} the classical
fact that ho\-mo\-lo\-gi\-cal dimension is additive on fibrations.
Theorem \ref{thm:top_gap} shows that a gap of length $g$ in the Betti
numbers of the fibre space implies that the total space has cohomology
in a dimension bigger than or equal to $g+1$.

\begin{center}
\rule{6ex}{0.1ex}
\end{center}

This paper supersedes the manuscript ``Homological identities for
Differential Gra\-ded Algebras, II'' from the spring of 2002.  That
manuscript suffers from technical problems which remain unsolved, and
it was never submitted.

Anders J.\ Frankild, my coauthor and friend of many years, died in
June 2007, before the present, more successful approach took the form
of this paper.

He will be bitterly missed, but his memory will live on.

Since Anders has not been able to check the final version of the
paper, the responsibility for any mistakes rests with me.

\section{Background}
\label{sec:background}

This section gives notation and e\-le\-men\-ta\-ry properties for DG
modules over simply connected cochain DG algebras.  For introductions
to the theory of DG algebras, see  \cite{AFH}, \cite{FHTbook}, or
\cite{KellerDG}.  The notation will stay close to \cite{FJiic},
\cite{PJampdgm}, and \cite{PJARtop}.

\begin{Setup}
\label{set:blanket}
By $k$ is denoted a field and by $R$ a cochain DG algebra over $k$
which has the form
\[
  \cdots \rightarrow 0 \rightarrow k \rightarrow 0 \rightarrow R^2
  \rightarrow R^3 \rightarrow \cdots
\]
and satisfies $\dim_k \H^i(R) < \infty$ for each $i$.
\end{Setup}

\begin{Remark}
In particular, $R^0 = k$, $R^1 = 0$, $\H^0(R) = k$ and $\H^1(R) = 0$. 
\end{Remark}

\begin{Notation}
There are derived categories $\D(R)$ of DG left-$R$-mo\-du\-les and
$\D(R^{\opp})$ of DG right-$R$-modules which support the derived
functors $\LTensor_R$ and $\RHom_R$.  We define full subcategories
of $\D(R)$,
\begin{align*}
  \Dp(R) & = \{\, M \in \D(R) \,|\, \mbox{$\H^j(M) = 0$ for $j \ll 0$} \,\}, \\
  \Dm(R) & = \{\, M \in \D(R) \,|\, \mbox{$\H^j(M) = 0$ for $j \gg 0$} \,\},
\end{align*}
and similarly for $\D(R^{\opp})$.

A DG left-$R$-module $M$ is {\em compact} precisely if it is finitely
built from $R$ in $\D(R)$ using distinguished triangles,
(de)suspensions, finite direct sums, and direct summands, cf.\ 
\cite[thm.\ 5.3]{KellerDG}.  The full subcategory of $\D(R)$
consisting of compact DG modules is denoted $\Dc(R)$, and similarly
for $\D(R^{\opp})$.

The suspension functor on DG modules is denoted by $\Sigma$.

The operation $(-)^{\natural}$ forgets the differential of a complex;
it sends DG algebras and DG modules to graded algebras and graded
modules. 

A DG $R$-module $M$ is called {\em locally finite} if it satisfies
$\dim_k \H^j(M) < \infty$ for each $j$.

The {\em infimum} and the {\em supremum} of a DG module are defined by
\[
  \inf M = \inf \{\, j \,|\, \H^j(M) \not= 0 \,\}, \;\;\;
  \sup M = \sup \{\, j \,|\, \H^j(M) \not= 0 \,\},
\]
and the {\em amplitude} is
\[
  \amp M = \sup M - \inf M.
\]
Note that we use the convention $\inf(\varnothing) = \infty$ and
$\sup(\varnothing) = -\infty$, so $\inf(0) = \infty$, $\sup(0) =
-\infty$, and $\amp(0) = -\infty$.  In fact, these special values
occur precisely when a DG module has zero cohomology, and this
property characterizes DG modules which are zero in the derived
category, so 
\begin{align}
\label{equ:infsupinfty}
\nonumber
  & \mbox{$M \cong 0$ in the derived category} \\
  & \qquad \Leftrightarrow \; \inf M = \infty 
  \; \Leftrightarrow \; \sup M = -\infty
  \; \Leftrightarrow \; \amp M = -\infty.
\end{align}
Moreover,
\begin{equation}
\label{equ:infleqsup}
  \mbox{$M \not\cong 0$ in the derived category}
  \; \Rightarrow \; \inf M \leq \sup M.
\end{equation}

We can view $k$ as a DG bi-$R$-module concentrated in cohomological
degree $0$.  The {\em Betti numbers} of a DG left-$R$-module $M$ are
\[
  \beta_R^j(M) = \dim_k \H^j(k \LTensor_R M).
\]

The {\em projective codimension} of $M$ is
\begin{equation}
\label{equ:pcd}
  \pcd_R M = \sup(k \LTensor_R M)
           = \sup \{\, j \,|\, \beta_R^j(M) \not= 0 \,\};
\end{equation}
it is an integer or $\infty$ or $-\infty$.
%Observe also that
%\begin{equation}
%\label{equ:l}
%  \inf(k \LTensor_R M) = \inf \{\, j \,|\, \beta^j(M) \not= 0 \,\};
%\end{equation}

If $\beta$ is a cardinal number, then a direct sum of $\beta$ copies
of $M$ will be denoted by $M^{(\beta)}$.
\end{Notation}

The following lemma holds by \cite{AFHpaper}.

\begin{Lemma}
\label{lem:semi-free}
Let $M$ be in $\Dp(R)$.  There is a semi-free resolution
\[
  \varphi : F \rightarrow M
\]
with semi-free filtration
\[
  0 = F(-1) \subseteq F(0) \subseteq F(1) \subseteq \cdots \subseteq F
\]
where the free quotients $F(j)/F(j-1)$ are direct sums of
(de)sus\-pen\-si\-ons $\Sigma^{\ell}R$ with $\ell \leq -\inf M$.
\end{Lemma}

%ALTERNATIVE to giving the reference:
%\begin{proof}
%This can be proved by combining the standard construction of semi-free
%resolutions, as given for instance in \cite[prop.\ 6.6(i)]{FHTbook},
%with $R^0 = k$ and $R^1 = 0$. 
%\end{proof}

\begin{Lemma}
\label{lem:inf2}
Let $P$ be in $\Dp(R^{\opp})$ and let $M$ be in $\Dp(R)$.  Then
\[
  \inf(P \LTensor_R M) = \inf P + \inf M.
\]
\end{Lemma}

\begin{proof}
If $P$ or $M$ is zero then the equation reads $\infty = \infty$, so
suppose that $P$ and $M$ are non-zero in the derived categories.  Then
$j = \inf P$ and $i = \inf M$ are integers.

By \cite[lem.\ 3.4(i)]{PJARtop}, we can replace $P$ with a
quasi-isomorphic DG mo\-du\-le which is zero in cohomological degrees
$<j$.  Lemma \ref{lem:semi-free} says that $M$ has a semi-free
resolution $F$ with a semi-free filtration where the successive
quotients are direct sums of DG modules $\Sigma^{\ell} R$ with $\ell
\leq -i$.  This implies that $F^{\natural}$ is a direct sum of graded
modules $\Sigma^{\ell} R^{\natural}$ with $\ell \leq -i$, so
$(P \otimes_R F)^{\natural}
= P^{\natural} \otimes_{R^{\natural}} F^{\natural}$
is a direct sum of graded modules $\Sigma^{\ell} P^{\natural}$ with
$\ell \leq -i$.  Since $P$ is zero in cohomological degrees $<j$, this
implies that $(P \otimes_R F)^{\natural}$ is zero in cohomological
degrees $<j+i$.  In particular we have $\inf(P \otimes_R F) \geq j+i$,
that is,
\begin{equation}
\label{equ:c}
  \inf(P \LTensor_R M) \geq j+i = \inf P + \inf M.
\end{equation}

On the other hand, a morphism of DG left-$R$-modules $\Sigma^{-i}R
\rightarrow M$ is determined by the image $z$ of $\Sigma^{-i}1_R$, and
$z$ is a cycle in $M^i$, the $i$th component of $M$.  Since
$\H^i(\Sigma^{-i}R) \cong \H^0(R) \cong k$, the induced map
$\H^i(\Sigma^{-i}R) \rightarrow \H^i(M)$ is just the map $k
\rightarrow \H^i(M)$ which sends $1_k$ to the cohomology class of $z$.
Hence if we pick cycles $z_{\alpha}$ such that the corresponding
cohomology classes form a $k$-basis of $\H^i(M)$ and construct a
morphism $\Sigma^{-i}R^{(\beta)} \rightarrow M$ by sending the
elements $\Sigma^{-i}1_R$ to the $z_{\alpha}$, then the induced map
$\H^i(\Sigma^{-i}R^{(\beta)}) \rightarrow \H^i(M)$ is an isomorphism.
Complete to a distinguished triangle
\begin{equation}
\label{equ:d}
  \Sigma^{-i}R^{(\beta)} \rightarrow M \rightarrow M^{\prime\prime} \rightarrow;
\end{equation}
since we have $\H^{i+1}(\Sigma^{-i}R^{(\beta)}) \cong
\H^1(R^{(\beta)}) = 0$, the long exact co\-ho\-mo\-lo\-gy sequence
shows 
\begin{equation}
\label{equ:k}
  \inf M^{\prime\prime} \geq i+1.
\end{equation}

Tensoring the distinguished triangle \eqref{equ:d} with $P$ gives
\[
  \Sigma^{-i}P^{(\beta)} \rightarrow P \LTensor_R M
    \rightarrow P \LTensor_R M^{\prime\prime} \rightarrow
\]
whose long exact cohomology sequence contains
\begin{equation}
\label{equ:e}
  \H^{j+i-1}(P \LTensor_R M^{\prime\prime})
    \rightarrow \H^{j+i}(\Sigma^{-i}P^{(\beta)})
    \rightarrow \H^{j+i}(P \LTensor_R M).
\end{equation}
The inequality \eqref{equ:c} can be applied to $P$ and
$M^{\prime\prime}$; because of \eqref{equ:k}, this gives
$\inf(P \LTensor_R M^{\prime\prime}) \geq j+i+1$
so the first term of the exact sequence \eqref{equ:e} is zero.  The
second term is 
$\H^{j+i}(\Sigma^{-i}P^{(\beta)}) \cong \H^j(P^{(\beta)})$
which is non-zero since $j = \inf P$.  So the third term is non-zero 
whence
\[
  \inf(P \LTensor_R M) \leq j+i = \inf P + \inf M.
\]
Combining with \eqref{equ:c} completes the proof.
\end{proof}

\begin{Lemma}
\label{lem:inf}
Let $P$ be non-zero in $\Dp(R^{\opp})$ and let $M$ be in $\Dp(R)$.
\begin{enumerate}

  \item  $M \cong 0$ in $\Dp(R)$ $\Leftrightarrow$ $P \LTensor_R M
    \cong 0$ in $\D(k)$.

\smallskip

  \item  $M \not\cong 0$ in $\Dp(R)$ $\Rightarrow$ $\inf M \leq \pcd M$.

\end{enumerate}
\end{Lemma}

\begin{proof}
(i)  Follows from Lemma \ref{lem:inf2} and Equation
\eqref{equ:infsupinfty}. 

(ii) When $M$ is non-zero in $\Dp(R)$, it follows from (i) that $k
\LTensor_R M$ is non-zero in $\D(k)$.  Then $\inf(k \LTensor_R M) \leq
\sup(k \LTensor_R M)$ by Equation \eqref{equ:infleqsup}.  By Lemma
\ref{lem:inf2} and Equation \eqref{equ:pcd}, this reads $\inf M \leq
\pcd M$.
\end{proof}

\section{A construction}
\label{sec:construction}

This section starts to set up the DG module adaptation of Serre's
ideas from the proof of \cite[thm.\ 10, p.\ 217]{Serre}.  The main
item is Construction \ref{con:basic} which approximates a DG module by
(de)suspensions of the DG algebra $R$.

\begin{Lemma}
\label{lem:basic}
Let $M$ be in $\Dp(R)$ and let $i$ be an integer with $i \leq
\inf M$.  There is a distinguished triangle in $\D(R)$,
\[
  \Sigma^{-i}R^{(\beta)} \rightarrow M \rightarrow M^{\prime\prime} \rightarrow,
\]
which satisfies the following.
\begin{enumerate}

  \item  $\beta = \beta^i(M)$.

\smallskip

  \item  $\beta^j(M) = \beta^j(M^{\prime\prime})$ for each $j \geq i+1$. 

\smallskip

  \item  $\inf M^{\prime\prime} \geq i+1$; in particular,
    $M^{\prime\prime}$ is in $\Dp(R)$. 

\end{enumerate}
If $M$ is locally finite then $\beta = \beta^i(M) < \infty$ and
$M^{\prime\prime}$ is also locally finite.
\end{Lemma}

\begin{proof}
The distinguished triangle is just \eqref{equ:d} from the proof of
Lemma \ref{lem:inf2}, constructed by picking cycles $z_{\alpha}$ in
$M^i$ such that the corresponding cohomology classes form a
$k$-basis for $\H^i(M)$, defining $\Sigma^{-i}R^{(\beta)}
\rightarrow M$ by sending the elements $\Sigma^{-i}1_R$ to the
$z_{\alpha}$, and completing to a distinguished triangle.

Property (iii) is the inequality \eqref{equ:k} in the proof of Lemma
\ref{lem:inf2}.  Ten\-so\-ring the distinguished triangle with $k$,
property (ii) is immediate and property (i) follows by using Lemma
\ref{lem:inf2}.

If $M$ is locally finite, then there are only finitely many cycles
$z_{\alpha}$, so $\beta < \infty$.  The long exact cohomology sequence
then shows that $M^{\prime\prime}$ is also locally finite.
\end{proof}

\begin{Construction}
\label{con:basic}
Let $M$ be non-zero in $\Dp(R)$ and write $i = \inf M$.  Observe that
$i$ is an integer.  Set $M \langle i \rangle = M$ and let $u \geq i$
be an integer.  By iterating Lemma \ref{lem:basic}, we can construct a
sequence of distinguished triangles in $\D(R)$,
\renewcommand{\arraystretch}{1.25}
\[
  \begin{array}{llllll}
  \Sigma^{-i}R^{(\beta^i)} & \rightarrow & M \langle i \rangle & \rightarrow & M \langle i+1 \rangle & \rightarrow, \\
  \Sigma^{-i-1}R^{(\beta^{i+1})} & \rightarrow & M \langle i+1 \rangle & \rightarrow & M \langle i+2 \rangle & \rightarrow, \\
  & & \;\;\;\;\;\;\; \vdots & & & \\
  \Sigma^{-u+1}R^{(\beta^{u-1})} & \rightarrow & M \langle u-1 \rangle & \rightarrow & M \langle u \rangle & \rightarrow, \\
  \Sigma^{-u}R^{(\beta^u)} & \rightarrow & M \langle u \rangle & \rightarrow & M \langle u+1 \rangle & \rightarrow, \\
  \end{array}
\]
\renewcommand{\arraystretch}{1}
where
\begin{enumerate}

  \item  $\beta^j = \beta^j(M)$ for each $j$.

\smallskip

  \item  $\beta^j(M) = \beta^j(M \langle \ell \rangle)$ for each $j \geq \ell$.

\smallskip

  \item  $\inf M \langle \ell \rangle \geq \ell$ for each $\ell$.  In
    particular, each $M \langle \ell \rangle$ is in $\Dp(R)$.

\end{enumerate}
If $M$ is locally finite, then so is each $M \langle \ell \rangle$, and
then each $\beta^j = \beta^j(M) = \beta^j(M \langle j \rangle)$ is
finite.
\end{Construction}

\begin{Proposition}
\label{pro:pcd}
\begin{enumerate}
  \item  If $M$ is in $\Dc(R)$, then it is locally finite and
    belongs to $\Dp(R)$.

\smallskip

  \item  Let $M$ be locally finite in $\Dp(R)$.  Then
\[
  \mbox{$M$ is in $\Dc(R)$ $\Leftrightarrow$ $\pcd M < \infty$.}
\]
\end{enumerate}
\end{Proposition}

\begin{proof}
Let $M$ be in $\Dc(R)$, that is, $M$ is finitely built from $R$ in
$\D(R)$.  Since $R$ is locally finite and belongs to $\Dp(R)$, the
same holds for $M$.  This proves (i).

Moreover, we have $\sup(k \LTensor_R R) = \sup k = 0 < \infty$ so we
must also have $\sup(k \LTensor_R M) < \infty$, that is, $\pcd M <
\infty$.  This proves (ii), implication $\Rightarrow$.

(ii), implication $\Leftarrow$: If $M$ is zero then it is certainly in
$\Dc(R)$, so assume that $M$ is non-zero in $\Dp(R)$.

Then $\inf M$ is an integer and $\inf M \leq \pcd M$ by Lemma
\ref{lem:inf}(ii).  On the other hand, $\pcd M < \infty$, so
\[
  \mbox{$p = \pcd M$ is an integer}.
\]
But $\pcd M = \sup(k \LTensor_R M)$ so $\beta^j(M) = \dim_k \H^j(k
\LTensor_R M) = 0$ for $j \geq p+1$.  In Construction \ref{con:basic},
by part (ii) this implies
\[
  \beta^j(M \langle p+1 \rangle) = 0 \;\mbox{for}\; j \geq p+1.
\]
However, part (iii) of the construction says $\inf M \langle p+1
\rangle \geq p+1$, so
\[
  \inf(k \LTensor_R M \langle p+1 \rangle) = \inf M \langle p+1 \rangle \geq p+1
\]
by Lemma \ref{lem:inf2}, that is
\[
  \beta^j(M \langle p+1 \rangle) 
    = \dim_k \H^j(k \LTensor_R M \langle p+1 \rangle) 
    = 0 \;\mbox{for}\; j < p+1.
\]
Altogether, $\beta^j(M \langle p+1 \rangle) = 0$ for each $j$.  That
is, each cohomology group of $k \LTensor_R M \langle p+1 \rangle$ is
zero and so $k \LTensor_R M \langle p+1 \rangle$ is itself zero.
Lemma \ref{lem:inf}(i) hence gives
\[
  M \langle p+1 \rangle \cong 0.
\]

But now the distinguished triangles in Construction \ref{con:basic},
starting with $\Sigma^{-p}R^{(\beta_p)} \rightarrow M \langle p
\rangle \rightarrow M \langle p+1 \rangle \rightarrow$ and running
backwards to the first one, $\Sigma^{-i}R^{(\beta_i)} \rightarrow M
\rightarrow M \langle i+1 \rangle \rightarrow$, show that $M$ is
finitely built from $R$ since each $\beta^j$ is finite.  That is, $M$
is in $\Dc(R)$.
\end{proof}

\begin{Remark}
\label{rmk:compact}
Let $M$ be non-zero in $\Dc(R)$.  Proposition \ref{pro:pcd}
implies that $M$ is locally finite in $\Dp(R)$ with $\pcd M < \infty$.
The proof of the proposition actually gives a bit more: First,
\[
  \mbox{$p = \pcd M$ is an integer.}
\]
Secondly, in Construction \ref{con:basic}, we have
\begin{equation}
\label{equ:iso_1}
  M \langle p+1 \rangle \cong 0
\end{equation}
in $\D(R)$.

Combining this isomorphism with the distinguished triangle 
\[
  \Sigma^{-p}R^{(\beta_p)} \rightarrow M \langle p \rangle
  \rightarrow M \langle p+1 \rangle \rightarrow
\]
shows
\begin{equation}
\label{equ:iso_2}
  M \langle p \rangle \cong \Sigma^{-p}R^{(\beta_p(M))}
\end{equation}
in $\D(R)$.
\end{Remark}

\section{Inequalities}
\label{sec:inequalities}

This section continues to set up the DG module adaptation of Serre's
ideas from the proof of \cite[thm.\ 10, p.\ 217]{Serre}.  The main
items are Propositions \ref{pro:inequalities} and
\ref{pro:compact_inequalities} which use Construction \ref{con:basic}
to prove some termwise inequalities of power series.

\begin{Setup}
\label{set:f}
Let
\[
  F : \D(R) \rightarrow \Mod(k)
\]
be a $k$-linear homological functor which respects coproducts.  For
$M$ in $\D(R)$, we set 
\[
  f_M(t) = \sum_{\ell} \dim_k F(\Sigma^{\ell} M)t^{\ell}.
\]
\end{Setup}

\begin{Remark}
In the generality of Setup \ref{set:f}, the expression $f_M(t)$ may
not belong to any reasonable set.  Some of the coefficients may be
infinite, and there may be non-zero coefficients in arbitrarily high
positive and negative degrees at the same time.

However, we will see that there are circumstances in which $f_M(t)$ is
a Laurent series.
\end{Remark}

\begin{Notation}
Termwise inequalities $\leq$ of coefficients between expressions
like $f_M(t)$ make sense; they will be denoted by $\preccurlyeq$.

Likewise, it makes sense to add such expressions, and to multiply them
by a number or a power of $t$.

Finally, the {\em degree} of an expression like $f_M(t)$ is defined
by 
\[
  \deg \Bigl( \sum_{\ell} f_{\ell}t^{\ell} \Bigr) = \sup \{\, \ell \,|\, f_{\ell} \not= 0 \,\}.
\]
\end{Notation}

\begin{Lemma}
\label{lem:inequality}
Let $M$ be in $\D(R)$.
\begin{enumerate}

  \item  $f_{\Sigma^j M}(t) = t^{-j}f_M(t)$.

\smallskip

  \item  $f_{M^{(\beta)}}(t) = \beta f_M(t)$.

\smallskip

  \item If
\[
  M^{\prime} \rightarrow M \rightarrow M^{\prime\prime} \rightarrow
\]
is a distinguished triangle in $\D(R)$, then there is a termwise
inequality
\[
  f_M(t) \preccurlyeq f_{M^{\prime}}(t) + f_{M^{\prime\prime}}(t).
\]
\end{enumerate}
\end{Lemma}

\begin{proof}
Parts (i) and (ii) are clear.  In part (iii), the distinguished
triangle gives a long exact sequence consisting of pieces
\[
  F(\Sigma^{\ell} M^{\prime}) 
  \rightarrow F(\Sigma^{\ell} M)
  \rightarrow F(\Sigma^{\ell} M^{\prime\prime}),
\]
whence
\[
  \dim_k F(\Sigma^{\ell} M) \leq
  \dim_k F(\Sigma^{\ell} M^{\prime}) + \dim_k F(\Sigma^{\ell} M^{\prime\prime})
\]
and the lemma follows.
\end{proof}

\begin{Proposition}
\label{pro:inequalities}
Let $M$ be non-zero in $\Dp(R)$.  Write $i = \inf M$, let $u \geq i$
be an integer, and consider Construction \ref{con:basic}.  There are
termwise inequalities
\begin{enumerate}
  \item  $f_M(t) \preccurlyeq
          (\beta^i(M)t^i + \cdots + \beta^u(M)t^u)f_R(t) + f_{M \langle u+1 \rangle}(t)$,

%\smallskip
%
%  \item  $\beta^i(M)t^if_R(t) \preccurlyeq
%          f_M(t) + t(\beta^{i+1}(M)t^{i+1} + \cdots +
%          \beta^u(M)t^u)f_R(t) + tf_{M \langle u+1 \rangle}(t)$,

\smallskip

  \item  $f_{M \langle u+1 \rangle}(t) \preccurlyeq
          f_M(t) + t^{-1}(\beta^i(M)t^i + \cdots + \beta^u(M)t^u)f_R(t)$.
\end{enumerate}
\end{Proposition}

\begin{proof}
This follows by applying Lemma \ref{lem:inequality} successively to
the distinguished triangles of Construction \ref{con:basic}.  For
instance, (i) can be proved as follows,
\begin{align*}
  f_M(t) & = f_{M \langle i \rangle}(t) \\
  & \preccurlyeq f_{\Sigma^{-i}R^{(\beta^i)}}(t) + f_{M \langle i+1 \rangle}(t) \\
  & = \beta^i(M)t^if_R(t) + f_{M \langle i+1 \rangle}(t) \\
  & \preccurlyeq \beta^i(M)t^if_R(t) + f_{\Sigma^{-i-1}R^{(\beta^{i+1})}}(t) + f_{M \langle i+2 \rangle}(t) \\
  & = \beta^i(M)t^if_R(t) + \beta^{i+1}(M)t^{i+1}f_R(t) + f_{M \langle i+2 \rangle}(t) \\
  & \preccurlyeq \cdots \\
  & = (\beta^i(M)t^i + \cdots + \beta^u(M)t^u)f_R(t) + f_{M \langle u+1 \rangle}(t),
\end{align*}
and (ii) is proved by similar manipulations.
\end{proof}

\begin{Proposition}
\label{pro:compact_inequalities}
Let $M$ be non-zero in $\Dc(R)$ and write $i = \inf M$ and
$p = \pcd M$.

Then $i$ and $p$ are integers with $i \leq p$, we have $\beta^p(M)
\not= 0$, and there are termwise inequalities
\begin{enumerate}
  \item  $f_M(t) \preccurlyeq
          (\beta^i(M)t^i + \cdots + \beta^p(M)t^p)f_R(t)$,

%\smallskip
%
%  \item  $\beta^i(M)t^if_R(t) \preccurlyeq
%          f_M(t) + t(\beta^{i+1}(M)t^{i+1} + \cdots + \beta^p(M)t^p)f_R(t)$,

\smallskip

  \item  $\beta^p(M)t^pf_R(t) \preccurlyeq
          f_M(t) + t^{-1}(\beta^i(M)t^i + \cdots + \beta^{p-1}(M)t^{p-1})f_R(t)$.
\end{enumerate}
\end{Proposition}

\begin{proof}
Proposition \ref{pro:pcd}(i) says that $M$ is in $\Dp(R)$, and
since $M$ is non-zero it follows that $i = \inf M$ is an integer.
Remark \ref{rmk:compact} says that $p = \pcd M$ is an integer.
Lemma \ref{lem:inf}(ii) says $i \leq p$.  Since $p = \sup(k \LTensor_R
M)$, it is clear that $\beta^p(M) = \dim_k \H^p(k \LTensor_R M) \not=
0$. 

Consider Construction \ref{con:basic} for $M$.  By Remark
\ref{rmk:compact}, Equations \eqref{equ:iso_1} and \eqref{equ:iso_2},
we have $M \langle p \rangle \cong \Sigma^{-p}R^{(\beta_p(M))}$ and $M
\langle p+1 \rangle \cong 0$.  Inserting this into the inequalities of
Proposition \ref{pro:inequalities} gives the inequalities of the
present proposition.
\end{proof}

As an immediate application, consider the following lemma.

\begin{Lemma}
\label{lem:identities}
Let $M$ be in $\Dc(R)$.  If $f_R(t)$ is a Laurent series in
$t^{-1}$ then so is $f_M(t)$, and
\[
  \deg f_M(t) = \deg f_R(t) + \pcd M.
\]
\end{Lemma}

\begin{proof}
If $M$ is zero then $f_M(t) = 0$ is trivially a Laurent series in
$t^{-1}$, and the equation of the lemma reads $-\infty = -\infty$ so
the lemma holds.

Suppose that $M$ is non-zero in $\Dc(R)$.  Since $f_R(t)$ is a Laurent
series in $t^{-1}$, Proposition \ref{pro:compact_inequalities}(i)
implies that so is $f_M(t)$ since each $\beta^j(M)$ is finite, cf.\ 
Proposition \ref{pro:pcd}(i) and Construction \ref{con:basic}.

If $f_R(t)$ has all coefficients equal to zero then Proposition
\ref{pro:compact_inequalities}(i) forces $f_M(t)$ to have all
coefficients equal to zero, and the equation of the lemma reads
$-\infty = -\infty$ so the lemma holds.

Suppose that not all coefficients of $f_R(t)$ are equal to zero.
Then Proposition \ref{pro:compact_inequalities}(i) implies
\[
  \deg f_M(t) \leq \deg f_R(t) + p = \deg f_R(t) + \pcd M.
\]
On the other hand, consider the inequality of Proposition
\ref{pro:compact_inequalities}(ii).  The left hand side contains a
non-zero monomial of degree $\deg f_R(t) + p$.  The right hand side
consists of two terms, and the second one, $t^{-1}(\beta^i(M)t^i +
\cdots + \beta^{p-1}(M)t^{p-1})f_R(t)$, consists of monomials of
degree $< \deg f_R(t)+ p$.  Hence the first term, $f_M(t)$, must
contain a non-zero monomial of degree $\deg f_R(t) + p$ whence
\[
  \deg f_M(t) \geq \deg f_R(t) + p = \deg f_R(t) + \pcd M.
\]
Combining the displayed inequalities gives the desired equation.
\end{proof}

\section{Main results}
\label{sec:identities}

This section shows the cochain Amplitude Inequality,
Aus\-lan\-der-Buchs\-baum Equality, and Gap Theorem for Bass numbers
in Corollaries \ref{cor:amp}, \ref{cor:AB}, and \ref{cor:Bass_gap}.
These results are special cases of Theorems \ref{thm:amp},
\ref{thm:AB}, and \ref{thm:gap}.

\begin{Setup}
From now on, we will only consider a special form of $F$ and $f$ from
Setup \ref{set:f}.  Namely, let $P$ be in $\D(R^{\opp})$ and set
\[
  F(-) = \H^0(P \LTensor_R -).
\]
This means that
\begin{equation}
\label{equ:fN}
  f_M(t) = \sum_{\ell} \dim_k \H^{\ell}(P \LTensor_R M)t^{\ell}
\end{equation}
and in particular
\begin{equation}
\label{equ:fR}
  f_R(t) = \sum_{\ell} \dim_k \H^{\ell}(P)t^{\ell}.
\end{equation}
\end{Setup}

\begin{Lemma}
\label{lem:AB}
Let $M$ be in $\Dc(R)$, and let $P$ be locally finite in
$\Dm(R^{\opp})$. 

Then
\[
  \sup(P \LTensor_R M) = \sup P + \pcd M.
\]
\end{Lemma}

\begin{proof}
The expression $f_R(t)$ is given by Equation \eqref{equ:fR} and it
is a Laurent series in $t^{-1}$ because $P$ is locally finite in
$\Dm(R^{\opp})$.  Lemma \ref{lem:identities} gives
\begin{equation}
\label{equ:a}
  \deg f_M(t) = \deg f_R(t) + \pcd M.
\end{equation}

However, Equations \eqref{equ:fN} and \eqref{equ:fR} imply
\[
  \deg f_M(t) = \sup(P \LTensor_R M)
\]
and
\[
  \deg f_R(t) = \sup P,
\]
so Equation \eqref{equ:a} reads
\[
  \sup(P \LTensor_R M) = \sup P + \pcd M
\]
as claimed.
\end{proof}

\begin{Theorem}
\label{thm:amp}
Let $M$ be in $\Dc(R)$.  Let $P$ be locally finite in $\D(R^{\opp})$
and suppose $\amp P < \infty$.  Then
\[
  \amp(P \LTensor_R M) = \amp P + \pcd M - \inf M.
\]
\end{Theorem}

\begin{proof}
Lemma \ref{lem:AB} says
\[
  \sup(P \LTensor_R M) = \sup P + \pcd M.
\]
Subtracting the equation of Lemma \ref{lem:inf2} produces
\[
  \sup(P \LTensor_R M) - \inf(P \LTensor_R M)
  = \sup P + \pcd M - \inf P - \inf M,
\]
and this is the equation of the present theorem.
\end{proof}

\begin{Corollary}
[Amplitude Inequality]
\label{cor:amp}
Let $M$ be non-zero in $\Dc(R)$.  Let $P$ be locally finite in
$\D(R^{\opp})$ and suppose $\amp P < \infty$.  Then
\[
  \amp(P \LTensor_R M) \geq \amp P.
\]
\end{Corollary}

\begin{proof}
Combine Theorem \ref{thm:amp} with Lemma \ref{lem:inf}(ii).
\end{proof}

\begin{Theorem}
\label{thm:AB}
Assume that there is a $P$ which is non-zero in $\Dc(R^{\opp})$ and
satisfies $\sup P < \infty$.  Set
\[
  d = \pcd P - \sup P.
\]

If $M$ is in $\Dc(R)$ and satisfies $\sup M < \infty$, then
\[
  \pcd M = \sup M + d.
\]
\end{Theorem}

\begin{proof}
Proposition \ref{pro:pcd}(i) gives that $M$ and $P$ are locally
finite.  They are also in $\Dm$ because they have finite supremum.
Hence Lemma \ref{lem:AB} says
\[
  \sup(P \LTensor_R M) = \sup P + \pcd M,
\]
and Lemma \ref{lem:AB} with $M$ and $P$ interchanged says
\[
  \sup(P \LTensor_R M) = \sup M + \pcd P.
\]

The two right hand sides must be equal,
\[
  \sup P + \pcd M = \sup M + \pcd P,
\]
and rearranging terms proves the proposition.
\end{proof}

\begin{Question}
For which DG algebras $R$ does there exist a DG module like $P$?  For
such DG algebras, the invariant $d = \pcd P - \sup P$ appears to be
interesting, and it would be useful to find a formula expressing it
directly in terms of $R$.
\end{Question}

The following corollary considers two easy, special cases of Theorem
\ref{thm:AB} which can reasonably be termed Auslander-Buchsbaum
Equalities.

\begin{Corollary}
[Auslander-Buchsbaum Equalities]
\label{cor:AB}
Let $M$ be in $\Dc(R)$.
\begin{enumerate}
  \item  If $R$ has $\sup R < \infty$, then
\[
  \pcd M = \sup M - \sup R.
\]

\smallskip

  \item  If $\sup M < \infty$ and $k$ is in $\Dc(R)$, then
\[
  \pcd M = \sup M + \pcd k.
\]
\end{enumerate}
\end{Corollary}

\begin{proof}
Both parts follow from Theorem \ref{thm:AB}, by using $R$ and $k$
in place of $P$.  In part (i), note that when $M$ is in $\Dc(R)$, it
is finitely built from $R$, so $\sup R < \infty$ implies $\sup M <
\infty$. 
\end{proof}

\begin{Theorem}
\label{thm:gap}
Let $M$ be locally finite in $\Dp(R)$.  Let $P$ be locally finite and
non-zero in $\D(R^{\opp})$ and suppose $\amp P < \infty$.

Let $g \geq \amp P$.  If the Betti numbers of $M$ have a gap of length
$g$ in the sense that there is a $j$ such that
\renewcommand{\arraystretch}{1.2}
\[
  \beta^{\ell}(M) \;
  \left\{
    \begin{array}{ll}
      \neq 0 & \mbox{for $\ell = j$}, \\
      = 0    & \mbox{for $j+1 \leq \ell \leq j+g$}, \\
      \neq 0 & \mbox{for $\ell = j+g+1$}, \\
    \end{array}
  \right.
\]
\renewcommand{\arraystretch}{1}
then
\[
  \amp(P \LTensor_R M) \geq g+1.
\]
\end{Theorem}

\begin{proof}
Since $\beta^j(M) \neq 0$, it is clear that $M$ is non-zero in $\Dp(R)$.
By (de)\-su\-spen\-ding, we can suppose $\inf M = \inf P = 0$.  Write
\[
  s = \sup P;
\]
then $s = \amp P$ and we have the assumption
\[
  g \geq s.
\]

Lemma \ref{lem:inf2} gives $\inf(P \LTensor_R M) = \inf P + \inf M =
0$, so $\H^0(P \LTensor_R M) \not= 0$.  To show the lemma, we need to
prove $\H^{\ell}(P \LTensor_R M) \not= 0$ for some $\ell \geq g+1$,
so let us assume
\[
  \H^{\geq g+1}(P \LTensor_R M) = 0
\]
and show a contradiction.  Lemma \ref{lem:inf2} implies
$\beta^{\ell}(M) = \dim_k \H^{\ell}(k \LTensor_R M) = 0$ for $\ell <
0$, so the integer $j$ from the proposition satisfies
$j \geq 0.$
Hence in particular
\begin{equation}
\label{equ:f}
  \H^{\geq j+g+1}(P \LTensor_R M) = 0.
\end{equation}

Inserting \eqref{equ:fN} and \eqref{equ:fR} into the inequality of
Proposition \ref{pro:inequalities}(i) gives
\begin{align*}
  \lefteqn{\sum_{\ell} \dim_k \H^{\ell}(P \LTensor_R M)t^{\ell}} \;\;\;\;\;\;\;\;\;\;\;\;\;\;\; & \\
  & \preccurlyeq (\beta^0(M)t^0 + \cdots \beta^u(M)t^u)  \sum_{\ell} \dim_k \H^{\ell}(P)t^{\ell}\\
  & \;\;\;\;\; + \sum_{\ell} \dim_k \H^{\ell}(P \LTensor_R M \langle u+1 \rangle)t^{\ell} \\
\end{align*}
where $u \geq 0$ is an integer and $M \langle u+1 \rangle$ is defined
by Construction \ref{con:basic}.  We have $\beta^{\ell}(M) = 0$ for
$j+1 \leq \ell \leq j+g$ while $\sum_{\ell} \dim_k
\H^{\ell}(P)t^{\ell}$ has terms only of degree $0, \ldots, s$, so the
first term on the right hand side is zero in degree $\ell$ for $j+s+1
\leq \ell \leq j+g$.  And $\inf M \langle u+1 \rangle \geq u+1$ by
Construction \ref{con:basic}(iii) so Lemma \ref{lem:inf2} implies
$\inf(P \LTensor_R M \langle u+1 \rangle) \geq u+1$, so by picking $u$
large we can move the second term on the right hand side into large
degrees and thereby ignore it.  It follows that the left hand side is
also zero in degree $\ell$ for $j+s+1 \leq \ell \leq j+g$; that is,
\[
  \mbox{$\H^{\ell}(P \LTensor_R M) = 0$ for $j+s+1 \leq \ell \leq j+g$.}
\]
Combining with Equation \eqref{equ:f} shows
\begin{equation}
\label{equ:g}
  \H^{\geq j+s+1}(P \LTensor_R M) = 0.
\end{equation}

Now insert \eqref{equ:fN} and \eqref{equ:fR} into Proposition
\ref{pro:inequalities}(ii) with $u = j$,
\begin{align*}
  \lefteqn{\sum_{\ell} \dim_k \H^{\ell}(P \LTensor_R M \langle j+1 \rangle)t^{\ell}} \;\;\;\;\;\;\;\;\;\;\;\;\;\;\;& \\
  & \preccurlyeq \sum_{\ell} \dim_k \H^{\ell}(P \LTensor_R M)t^{\ell} \\
  & \;\;\;\;\; + t^{-1}(\beta^0(M)t^0 + \cdots + \beta^j(M)t^j)
               \sum_{\ell} \dim_k \H^{\ell}(P)t^{\ell}. \\
\end{align*}
Again, $\sum_{\ell} \dim_k \H^{\ell}(P)t^{\ell}$ only has terms of
degree $0, \ldots, s$, so on the right hand side, the second term is
zero in degrees $\geq j+s$.  In particular, it is zero in degrees
$\geq j+s+1$, and since Equation \eqref{equ:g} implies that the same holds
for the first term, it must also hold for the left hand side, that is,
\begin{equation}
\label{equ:h}
  \H^{\geq j+s+1}(P \LTensor_R M \langle j+1 \rangle) = 0.
\end{equation}

Now, $\inf M \langle j+1 \rangle \geq j+1$ by Construction
\ref{con:basic}(iii), so Lemma \ref{lem:inf2} implies
\[
  \inf(k \LTensor_R M \langle j+1 \rangle) \geq j+1.
\]
And Construction \ref{con:basic}(ii) says $\beta^{\ell}(M \langle
j+1 \rangle) = \beta^{\ell}(M)$ for $\ell \geq j+1$, so
$\beta^{\ell}(M) = 0$ for $j+1 \leq \ell \leq j+g$ gives
\[
  \mbox{$\beta^{\ell}(M \langle j+1 \rangle) = 0$ for $j+1 \leq \ell \leq j+g$},
\]
that is,
\[
  \mbox{$\H^{\ell}(k \LTensor_R M \langle j+1 \rangle) = 0$ for $j+1 \leq \ell \leq j+g$},
\]
so we even have
\[
  \inf(k \LTensor_R M \langle j+1 \rangle) \geq j+g+1,
\]
that is, $\inf M \langle j+1 \rangle \geq j+g+1$ by Lemma
\ref{lem:inf2} again.  Hence
\[
  \inf(P \LTensor_R M \langle j+1 \rangle) \geq j+g+1
\]
by Lemma \ref{lem:inf2}.  However, $g \geq s$, so the only way this
can be compatible with Equation \eqref{equ:h} is if we have
\[
  P \LTensor_R M \langle j+1 \rangle \cong 0.
\]
By Lemma \ref{lem:inf}(i) this means that $M \langle j+1 \rangle \cong
0$.  And this gives 
\[
  \beta^{j+g+1}(M) = \beta^{j+g+1}(M \langle j+1 \rangle) = 0,
\]
which is the desired contradiction since we had assumed
$\beta^{j+g+1}(M) \not= 0$.
\end{proof}

\begin{Corollary}
[Gap Theorem for Betti numbers]
\label{cor:Betti_gap}
Suppose $\sup R < \infty$.  Let $M$ be locally finite in $\Dp(R)$.

Let $g \geq \sup R$.  If the Betti numbers of $M$ have a gap of length
$g$ in the sense that there is a $j$ such that
\renewcommand{\arraystretch}{1.2}
\[
  \beta^{\ell}(M) \;
  \left\{
    \begin{array}{ll}
      \neq 0 & \mbox{for $\ell = j$}, \\
      = 0    & \mbox{for $j+1 \leq \ell \leq j+g$}, \\
      \neq 0 & \mbox{for $\ell = j+g+1$}, \\
    \end{array}
  \right.
\]
\renewcommand{\arraystretch}{1}
then
\[
  \amp M \geq g+1.
\]
\end{Corollary}

\begin{proof}
This follows from Theorem \ref{thm:gap} by using $R$ in
place of $P$.
\end{proof}

\begin{Remark}
Conversely, if $\sup R < \infty$ and $M$ is locally finite in $\Dp(R)$
with $\amp M \leq \sup R$, then the Betti numbers of $M$ can have no
gaps of length bigger than or equal to $\sup R$.
\end{Remark}

By evaluating the previous theorem on the $k$-linear dual
$\Hom_k(M,k)$, we immediately get the following result in which the
{\em Bass numbers} of a DG module are
\[
  \mu^j(M) = \dim_k \H^j(\RHom_R(k,M)).
\]

% \begin{Theorem}
% [Bass]
% \label{thm:Bass}
% Let $M$ in $\D(R)$ have $\Hom_k(M,k)$ in $\Dc(R^{\opp}) and $-\infty
% < \inf M$.  
% \begin{enumerate}
%   \item  If $R$ has $\sup R < \infty$, then
% \[
%   \depth M = \inf M + \sup R.
% \]

% \smallskip

%   \item  If $k$ is in $\Dc(R)$, then
% \[
%   \depth M = \inf M - \pcd k.
% \]
% \end{enumerate}
% \end{Theorem}

\begin{Corollary}
[Gap Theorem for Bass numbers]
\label{cor:Bass_gap}
Suppose $\sup R < \infty$.  Let $M$ be locally finite in $\Dm(R)$.

Let $g \geq \sup R$.  If the Bass numbers of $M$ have a gap of length
$g$ in the sense that there is a $j$ such that
\renewcommand{\arraystretch}{1.2}
\[
  \mu^{\ell}(M) \;
  \left\{
    \begin{array}{ll}
      \neq 0 & \mbox{for $\ell = j$}, \\
      = 0    & \mbox{for $j+1 \leq \ell \leq j+g$}, \\
      \neq 0 & \mbox{for $\ell = j+g+1$}, \\
    \end{array}
  \right.
\]
\renewcommand{\arraystretch}{1}
then
\[
  \amp M \geq g+1.
\]
\end{Corollary}

\begin{Remark}
\label{rmk:Bass_gap}
Conversely, if $\sup R < \infty$ and $M$ is locally finite in
$\Dm(R)$ with $\amp M \leq \sup R$, then the Bass numbers
of $M$ can have no gaps of length bigger than or equal to $\sup R$.

In particular, the Bass numbers of $R$ itself can have no gaps of
length bigger than or equal to $\sup R$.  This shows for the present
class of DG algebras that the answer is affirmative to the question
asked by Avramov and Foxby in \cite[Question 3.10]{AF} for local chain
DG algebras.  In fact, it shows that for simply connected cochain DG
algebras, Avramov and Foxby's conjectural bound on the gap length of
the Bass numbers can be sharpened by an amount of one.
\end{Remark}

\section{Topology}
\label{sec:topology}

This section applies the Auslander-Buchsbaum Equality and the Gap
Theorem to the singular cochain DG algebra of a topological space.
The context will be a fibration of topological spaces, and we recover
in Theorem \ref{thm:additivity} that ho\-mo\-lo\-gi\-cal dimension is
additive on fibrations.  Theorem \ref{thm:top_gap} shows that a gap of
length $g$ in the Betti numbers of the fibre space implies that the
total space has non-zero cohomology in a dimension $\geq g+1$.

A reference for the algebraic topology of this section is
\cite{FHTbook}.

\begin{Setup}
Let
\[
  F \rightarrow X \rightarrow Y
\]
be a fibration of topological spaces where $\dim_k \H^j(X;k) < \infty$
and $\dim_k \H^j(Y;k) < \infty$ for each $j$ and where $Y$ is simply
connected.
\end{Setup}

\begin{Remark}
Recall that the {\em singular cohomology} $\H^j(Z;k)$ of a
topological space $Z$ is defined in terms of the {\em singular cochain
complex} $\C^*(Z;k)$ by
\[
  \H^j(Z;k) = \H^j(\C^*(Z;k)).
\]

The singular cochain complex is a DG algebra, and by \cite[exa.\ 6,
p.\ 146]{FHTbook} the assumptions on the space $Y$ mean that
$\C^*(Y;k)$ is quasi-isomorphic to a DG algebra which falls under
Setup \ref{set:blanket}, so the results proved so far apply to it.

Moreover, the continuous map $X \rightarrow Y$ induces a morphism
$\C^*(Y;k) \rightarrow \C^*(X;k)$ whereby $\C^*(X;k)$ becomes a DG
bi-$\C^*(Y;k)$-module which is locally finite and belongs to $\Dp$ by
the assumptions on $X$.
\end{Remark}

\begin{Notation}
The dimensions $\dim_k \H^j(Z;k)$ are called the {\em Betti numbers}
of the topological space $Z$.

By
\begin{equation}
\label{equ:hd}
  \hd Z = \sup \{\, j \,|\, \H^j(Z;k) \not= 0 \,\} = \sup \C^*(Z;k)
\end{equation}
is denoted the {\em homological dimension} of $Z$; it is a
non-negative integer or $\infty$.

By $\Omega Z$ is denoted the {\em Moore loop space} of $Z$.
\end{Notation}

\begin{Lemma}
\label{lem:hd}
\begin{enumerate}

  \item  We have
\[
  \dim_k \H^j(F;k) = \beta^j_{\C^*(Y;k)}(\C^*(X;k))
\]
and $\hd F = \pcd_{\C^*(Y;k)}(\C^*(X;k))$.

\smallskip

  \item  We have
\[
  \dim_k \H^j(\Omega Y;k) = \beta^j_{\C^*(Y;k)}(k)
\]
and $\hd \Omega Y = \pcd_{\C^*(Y;k)}(k)$.

\end{enumerate}
\end{Lemma}

\begin{proof}
(i) We know that $\C^*(F;k) \cong k \LTensor_{\C^*(Y;k)} \C^*(X;k)$ in
$\D(k)$ by \cite[thm.\ 7.5]{FHTbook}.  Taking the dimension of the
$j$'th cohomology proves the displayed equation, and the other
equation is an immediate consequence, cf.\ Equations \eqref{equ:pcd}
and \eqref{equ:hd}.

(ii) This follows by using (i) on the fibration $\Omega Y \rightarrow
PY \rightarrow Y$ where $PY$ is the {\em Moore path space} of $Y$: The
space $PY$ is contractible so $\C^*(PY;k)$ is isomorphic to the DG
module $k$ in $\D(\C^*(Y;k))$.
\end{proof}

\begin{Theorem}
[Additivity of homological dimension]
\label{thm:additivity}
\begin{enumerate}

  \item  If $\hd F < \infty$ and $\hd Y < \infty$, then
\[
  \hd X = \hd F + \hd Y.
\]

\smallskip

  \item  If $\hd X < \infty$ and $\hd \Omega Y < \infty$, then
\[
  \hd X = \hd F - \hd \Omega Y.
\]

\end{enumerate}
\end{Theorem}

\begin{proof}
Let us apply Corollary \ref{cor:AB} to the data
\[
  \mbox{$R = \C^*(Y;k)$ \;and\; $M = \C^*(X;k)$.}
\]

(i) Using Lemma \ref{lem:hd}(i) we have
\[
  \pcd_R(M) = \pcd_{\C^*(Y;k)}(\C^*(X;k)) = \hd F < \infty.
\]
By Proposition \ref{pro:pcd}(ii), this says that $M$ is in $\Dc(R)$.
Moreover, Equation \eqref{equ:hd} gives
\[
  \sup R = \sup \C^*(Y;k) = \hd Y < \infty.
\]

This shows that Corollary \ref{cor:AB}(i) does apply, and evaluating its
equation gives $\hd F = \hd X - \hd Y$, proving (i).

(ii) Using Lemma \ref{lem:hd}(ii) we have
\[
  \pcd_R(k) = \pcd_{\C^*(Y;k)}(k) = \hd \Omega Y < \infty.
\]
By Proposition \ref{pro:pcd}(ii), this says that $k$ is in $\Dc(R)$.
But we also have
\[
  \sup M = \sup \C^*(X;k) = \hd X < \infty,
\]
and since $M$ is locally finite, it follows that $\dim_k \H(M) <
\infty$ whence $M$ is finitely built from $k$ in $\D(R)$.  Hence $M$
is also in $\Dc(R)$.

This shows that Corollary \ref{cor:AB}(ii) does apply, and evaluating
its equation gives $\hd F = \hd X + \hd \Omega Y$, proving (ii).
\end{proof}

\begin{Theorem}
[Gap]
\label{thm:top_gap}
\begin{enumerate}

  \item  If $\hd Y < \infty$ and the Betti numbers of $F$ have a
    gap of length $g \geq \hd Y$ in the sense that there is a $j$
    such that
\renewcommand{\arraystretch}{1.2}
\[
  \H^{\ell}(F;k) \;
  \left\{
    \begin{array}{ll}
      \neq 0 & \mbox{for $\ell = j$}, \\
      = 0    & \mbox{for $j+1 \leq \ell \leq j+g$}, \\
      \neq 0 & \mbox{for $\ell = j+g+1$}, \\
    \end{array}
  \right.
\]
\renewcommand{\arraystretch}{1}
then
\[
  \hd X \geq g+1.
\]

\smallskip

  \item  If $\hd X < \infty$ and the Betti numbers of $\Omega Y$
    have a gap of length $g \geq \hd X$ in the sense that there is
    a $j$ such that
\renewcommand{\arraystretch}{1.2}
\[
  \H^{\ell}(\Omega Y;k) \;
  \left\{
    \begin{array}{ll}
      \neq 0 & \mbox{for $\ell = j$}, \\
      = 0    & \mbox{for $j+1 \leq \ell \leq j+g$}, \\
      \neq 0 & \mbox{for $\ell = j+g+1$}, \\
    \end{array}
  \right.
\]
\renewcommand{\arraystretch}{1}
then
\[
  \hd F \geq g+1,
\]
and $\hd F < \infty$ forces $\hd Y = \infty$.

\end{enumerate}
\end{Theorem}

\begin{proof}
(i) Let us apply Corollary \ref{cor:Betti_gap}(i) to the data
\[
  \mbox{$R = \C^*(Y;k)$ \;and\; $M = \C^*(X;k)$.}
\]
Then $\sup R = \hd Y$ is clear, $\beta^j(M) = \dim_k \H^j(F;k)$
holds by Lemma \ref{lem:hd}(i), and $\amp M = \hd X$ is clear, so (i)
follows.

(ii) Let us apply Theorem \ref{thm:gap} to the data
\[
  \mbox{$R = \C^*(Y;k)$,\; $M = k$, \;and\; $P = \C^*(X;k)$.}
\]
Then $\amp P = \hd X$ is clear, $\beta^j(M) = \beta^j(k) = \dim_k
\H^j(\Omega Y;k)$ holds by Lemma \ref{lem:hd}(ii), and
\[
  \amp(P \LTensor_R M)
  = \amp(\C^*(X;k) \LTensor_{\C^*(Y;k)} k)
  = \amp(\C^*(F;k))
  = \hd F
\]
since $\C^*(X;k) \LTensor_{\C^*(Y;k)} k \cong \C^*(F;k)$ by
\cite[thm.\ 7.5]{FHTbook}, so the inequality of (ii) follows.

Moreover, if we had $\hd F < \infty$ and $\hd Y < \infty$, then
Theorem \ref{thm:additivity}(i) would apply, and we would get the
contradiction
\[
  g \geq \hd X = \hd F + \hd Y \geq g + 1 + \hd Y \geq g+1.
\]
\end{proof}

\begin{Example}
Set $Y = S^n \vee S^{n+1} \vee S^{n+2} \vee \cdots$ for an $n \geq
2$.  Then $Y$ is simply connected with $\dim_{\BQ} \H^j(Y;\BQ) <
\infty$ for each $j$, and
\[
  \sum_j \dim_{\BQ} \H^j(\Omega Y;\BQ)t^j = 
  \frac{1}{1 - (t^n + t^{n+1} + t^{n+2} + \cdots)}
\]
by \cite[exa.\ 1, p.\ 460]{FHTbook}.  It follows that $\H^1(\Omega
Y;\BQ) = \cdots = \H^{n-1}(\Omega Y;\BQ) = 0$, so the Betti numbers of
$\Omega Y$ have a gap of length $n-1$.

Hence, if $F \rightarrow X \rightarrow Y$ is a fibration with
$\dim_{\BQ} \H^*(X;\BQ) < \infty$ and $\hd X \leq n-1$, then Theorem
\ref{thm:top_gap}(ii) says $\hd F \geq n$.
\end{Example}

%\medskip
%\noindent
%{\em Acknowledgement. }

\end{document}